\documentclass[12pt]{amsart}
\usepackage{amssymb}
\usepackage{color}
\usepackage{amsmath}
\usepackage{mathtools}
\usepackage[all]{xy}
\usepackage{enumitem}
\usepackage{tabularx}
\usepackage{booktabs}

 \rm

\renewcommand{\(}{\left(}
\renewcommand{\)}{\right)}

\newcommand{\<}{\left\langle}
\renewcommand{\>}{\right\rangle}
\newcommand{\x}{\times}
\renewcommand{\bar}{\overline}
\newcommand{\abs}[1]{\left\lvert#1\right\rvert}
\newcommand{\norm}[1]{\left\lVert#1\right\rVert}
\newcommand{\inj}[1]{\operatorname{inj}_{#1}}
\newcommand{\grad}[1]{\operatorname{grad}_{#1}}
\newcommand{\st}{\:|\:}

\renewcommand{\phi}{\varphi}

\newcommand{\R}{{\mathbb{R}}}

\theoremstyle{plain}
\newtheorem{thm}{Theorem}[section]
\newtheorem{lem}[thm]{Lemma}

\theoremstyle{definition}

\theoremstyle{remark}

\title{A generalization of a result of Minakshisundaram and Pleijel}

\author{ Ankita~Sharma, Mansi~Mishra and M.~K.~Vemuri}
\thanks{The first author was supported by a fellowship from CSIR
(File No.: 09/1217(0077)/2019-EMR-I)}
\address{Department of Mathematical Sciences, IIT(BHU), Varanasi, 221005,
INDIA}

\begin{document}


\begin{abstract}
Minakshisundaram and Pleijel gave an asymptotic formula for
the sum of squares of the pointwise values of the eigenfunctions of
the Laplace-Beltrami operator on a compact Riemannian manifold, with
eigenvalues less than a fixed number.  Zelditch later extended this 
result by replacing the pointwise values with the Fourier coefficients of
a smooth measure supported on a compact submanifold.  Zelditch's result
is very general, and his proof relies on the theory of Fourier integral
operators.

Here we give a proof based on methods of Riemannian geometry.
\end{abstract}

\maketitle
\thispagestyle{empty}

\section{Introduction}\label{S:intro}
Let $M$ be a compact Riemannian manifold of dimension $m$ and let
$\mu$ denote the Riemannian measure on $M$.
Let $\Delta$ denote the {\em Laplace-Beltrami operator} of $M$.
It is well known (see e.g., \cite{roe}) that
there exist $\lambda_j \in [0,\infty)$, and $\phi_j \in C^\infty(M)$,
$j=0, 1, 2, \dots$, such that $\{\phi_j\}_{j=0}^\infty$ is an orthonormal
basis of $L^2(M, \mu)$ and
\begin{equation*}
-\Delta \phi_j = \lambda_j \phi_j.
\end{equation*}

Thus, by Parseval's formula (see e.g., \cite[Theorem 4.18]{Rudin}),
if $\psi\in C^\infty(M)$ then
\begin{equation}\label{E:parseval}
\sum_{\lambda_j < T}\abs{\widehat{\psi}(j)}^2 \sim 
\int_M \abs{\psi}^2 \,d\mu, \quad T\to\infty,
\end{equation}
where $\widehat{\psi}(j)=\<\psi,\phi_j\>$ is the $j$-th Fourier
coefficient of $\psi$.

On the other hand,
in \cite{minakshisundaram}, Minakshisundaram and Pleijel proved that
for all $P\in M$,
\begin{equation}\label{E:MP2}
\sum_{\lambda_j < T}\abs{\phi_j(P)}^2 \sim 
\frac{T^{m/2}}{(4\pi)^{m/2}\Gamma(\frac{m}{2}+1)},\quad T\to \infty.
\end{equation}
We may interpret this as a statement about the Fourier coefficients of
the delta measure at $P$.
Recall that if $\tau$ is a measure on $M$, the $j$-th Fourier coefficient
of $\tau$ (as a distribution on $M$) is
\begin{equation*}
\widehat{\tau}(j) = \<\tau, \phi_j\> = \int_M \phi_j \,d\tau.
\end{equation*}
In particular, if $\delta_P$ is the delta measure
at $P$, then
\begin{equation*}
\widehat{\delta_P}(j) = \int_M\phi_j \,d\delta_P
                      = \phi_j(P).
\end{equation*}
Thus (\ref{E:MP2}) may be restated as
\begin{equation*}
\sum_{\lambda_j < T} \abs{\widehat{\delta_P}(j)}^2 \sim 
\frac{T^{m/2}}{(4\pi)^{m/2}\Gamma(\frac{m}{2}+1)},\quad T\to \infty.
\end{equation*}

In this work, we show that there is an analogous asymptotic formula
for the Fourier coefficients of any smooth measure supported on
a submanifold of $M$, which has
both (\ref{E:parseval}) and (\ref{E:MP2}) as special cases.

Suppose $N$ is a compact submanifold of $M$ and let
$\nu$ denote the Riemannian measure on $N$.  Let $\psi \in C^\infty(N)$.
A measure of the form $\tau = \psi\nu$ is called a smooth measure supported
on $N$ (see e.g. \cite[Chapter 8, \S3]{stein}).  
The main result of this work is the following theorem.

\begin{thm}\label{T:main}
Let $\tau=\psi\nu$ be a smooth measure supported on a compact
codimension $k$ submanifold $N$ of $M$.
Then
\begin{equation*}
\sum_{\lambda_j < T}\abs{\widehat{\tau}(j)}^2 \sim 
\frac{T^{k/2}\int_N\abs{\psi}^2 d\nu}{(4\pi)^{k/2}\Gamma(\frac{k}{2}+1)},
\quad T\to \infty.
\end{equation*}
\end{thm}

In \cite[Corollary 3.3]{steven}, Zelditch gave a very general result
about Fourier integral operators, from which Theorem \ref{T:main} may
be deduced.  Here we give a proof of Theorem \ref{T:main} in the spirit of
Minakshisundaram and Pleijel \cite{minakshisundaram}.
We give an application of Theorem \ref{T:main} to an asymptotic Pythagorean
identity for the Legendre polynomials.
Analogous identities for other families of orthogonal polynomials will appear
in a forthcoming work.

\subsection{Acknowledgements}
We would like to thank C.\ S.\ Aravinda for helpful comments on this work.
The third author would like to thank Mokshay Madiman for discussions
related to this work during his numerous visits to the University of Delaware.
We would also like to thank the anonymous referee for a suggestion which
substantially simplified the original argument.

\section{Discussion of the main result}

Theorem \ref{T:main} follows from an asymptotic expansion at $0$ of
the $L^2$ norm of the heat flow of $\tau$.

Let $k_t(x,y)$ denote the heat kernel of $M$. The heat flow of $\tau$ is
given by
\begin{equation*}
f_t(x)
= \int_M k_t(x,y) \, d\tau(y) = \int_N k_t(x,y)\psi(y) \, d\nu(y),
\quad x\in M, t>0.
\end{equation*}
In Section \ref{S:heat-flow}, we prove that
\begin{equation}\label{E:heat-flow}
\norm{f_{t/2}}_2^2 \sim (4\pi t)^{-k/2}\norm{\psi}_2^2, \quad t\to 0.
\end{equation}

Since
\begin{equation*}
k_t(x,y)=\sum_{j=0}^\infty e^{-\lambda_j t} \phi_j(x)\phi_j(y),
\quad x,y\in M, t>0,
\end{equation*}
it follows that
\begin{equation*}
\norm{f_{t/2}}_2^2
=\sum_{j=0}^\infty \abs{\hat{\tau}(j)}^2 e^{-\lambda_j t}
=\int_0^\infty e^{-tT} \, d\alpha(T),
\end{equation*}
where
$\alpha:[0,\infty) \to \R$ is given by
\begin{equation*}
\alpha(T) = \sum_{\lambda_j<T}\abs{\hat{\tau}(j)}^2.
\end{equation*}
In other words, $\norm{f_{t/2}}_2^2$ is the Laplace transform of the measure
$d\alpha$.
Therefore Theorem \ref{T:main} follows from Equation (\ref{E:heat-flow}) and
the following result (see \cite[Chapter 5, Theorem 4.3]{widder}).

\begin{thm}[Karamata's Tauberian theorem]
If $\alpha(t)$ is non-decreasing and such that the integral
\begin{equation*}
f(s)=\int_0^\infty e^{-st}\,d\alpha(t)
\end{equation*}
converges for $s>0$, and if for some non-negative number $\gamma$
\begin{equation*}
f(s)\sim \frac{A}{s^\gamma}, \quad s\to 0+
\end{equation*}
then
\begin{equation*}
\alpha(t)\sim \frac{At^\gamma}{\Gamma(\gamma+1)}, \quad t\to \infty.
\end{equation*}
\end{thm}

\subsection{Outline of the paper}

In Section \ref{S:distance-function}, we derive a formula for the Hessian
on $N$ of the squared distance function to a point in $M$ in terms of the
second fundamental form of $N$ and that of an osculating geodesic sphere.
We use this formula to estimate the Hessian on $N$ of the squared
distance function.  These estimates allow us to use Laplace's method
to find the asymptotic of a Gaussian-type integral involving the
distance function.
In Section \ref{S:heat-flow}, we use the results of Section
\ref{S:distance-function}
to derive Equation (\ref{E:heat-flow}).
In Section \ref{S:lp}, we illustrate Theorem \ref{T:main} by deriving an
asymptotic Pythagorean identity for the Legendre polynomials.

\subsection{Notation}

For the convenience of the reader, we summarize the notation used
throughout this work.

\medskip

\begin{center}
\begin{tabularx}{0.75\textwidth}{p{0.16\textwidth}X}
\toprule
\bf{Symbol}       &  \bf{Description}\\
\toprule
$\overline\nabla$ &  Levi-Civita connection on $M$\\
$\mu$             &  Riemannian measure on $M$\\
$d_M$             &  Distance function on $M$\\
$\inj{M}$         &  Injectivity radius of $M$\\
$R$               &  Curvature tensor of $M$\\
$k_t$             &  Heat kernel of $M$\\
$\nabla$          &  Levi-Civita connection on $N$\\
$\nu$             &  Riemannian measure on $N$\\
$d_N$             &  Distance function on $N$\\
$\inj{N}$         &  Injectivity radius of $N$\\
$\exp$            &  Exponential map of $N$\\
$\vec{h}$         &  Second fundamental form of $N$ in $M$\\
$\boldsymbol{\nu}$&  Normal bundle of $N$ in $M$\\
$\exp_{\boldsymbol{\nu}}$& Exponential map of $\boldsymbol{\nu}$\\
$N_\eta$           &  $\{x \in M \st d_M(x,N)<\eta \}$\\
$d_M(x,N)$        &  Distance from $x$ to $N$\\
\bottomrule
\end{tabularx}
\end{center}

\section{The distance function on $N$ to a point of $M$}\label{S:distance-function}

Let $n$, $\nabla$, $\exp$, $d_N$ and $\inj{N}$ denote the
dimension,
Levi-Civita connection, exponential map, distance function and
injectivity radius of $N$ respectively.
We will use the same symbol $\nabla$ for the connection on all tensor
bundles over $N$; in particular, if $g$ is a function on $N$, $\nabla g$ denotes the
differential of $g$ on $N$.
Note that $\inj{N}>0$ because $N$ is compact.
Let $\vec{h}$ denote the vector-valued second fundamental form
of $N$ in $M$.
Since $N$ is compact, $\vec{h}$ is bounded, i.e., there exists a
constant $\kappa$ such that for all $y\in N$ and $X,Y\in T_yN$,
\begin{equation*}
\norm{\vec{h}(X,Y)} \leq \kappa \norm{X}\norm{Y}.
\end{equation*}

For $\eta\in (0,\infty)$, let $N_{\eta}$ denote the set
$\{x \in M \st d_M(x,N)<\eta \}$.
Observe that there exists $\eta\in (0, \inj{M})$ such that the set
$N_\eta$ is a tubular neighborhood of $N$, i.e.,
for every $x \in N_\eta$ there exists a unique $z \in N$ with
$d_M(x,z)=d_M(x,N)$,
and the map $\Pi:N_\eta \to N$ which sends $x$ to $z$ is smooth
\cite[Chapter 7, Proposition 26]{tubular}.
For a fixed $x \in N_\eta$, define $\rho_{_x} : M \to \R$ by
\begin{equation*}
\rho_{_x}(y) = {d_M(x,y)}^2-{d_M(x,N)}^2.
\end{equation*}

\begin{lem}\label{L:hessian}
For all $x\in N_\eta$ and unit vectors $Y\in T_{\Pi(x)}N$,
\begin{equation*}
\begin{gathered}
(\nabla^2\rho_{_x})_{\Pi(x)}(Y,Y) \le
\frac{2d_M(x,N)\sqrt{\lambda}}{\tanh(d_M(x,N)\sqrt{\lambda})}+ 2\kappa d_M(x,N),
\quad \text{and}\\
(\nabla^2\rho_{_x})_{\Pi(x)}(Y,Y) \ge
\frac{2d_M(x,N)\sqrt{\lambda}}{\tan(d_M(x,N)\sqrt{\lambda})}- 2\kappa d_M(x,N).
\end{gathered}
\end{equation*}

\end{lem}  

\begin{proof}
Let $x\in N_{\eta}$, $L=\{y \in M \:|\: d_M(x,y)= d_M(x,N)\}$, and $\tilde{h}$ the second fundamental
form of $L$ with respect to the outward pointing unit normal.
Observe that $L$ is smooth and $T_{\Pi(x)}N \subseteq T_{\Pi(x)}L$.  Let
$X,Y \in T_{\Pi(x)}N$.
Extend $X$ and $Y$ to be smooth vector fields on $M$ tangent to $L$.
Let $\bar{\nabla}$ denote the Levi-Civita connection on $M$.  Again,
the same symbol $\bar{\nabla}$ is used for the
connection on all tensor bundles over $M$.  Then
\begin{equation*}
\begin{gathered}
\nabla^2\rho_{_x}(X,Y)
 = X(Y\rho_{_x}) - (\nabla_X Y)\rho_{_x}, \quad \text{and}\\
\bar{\nabla}^2 \rho_{_x} (X,Y)
 = X(Y\rho_{_x}) - (\bar{\nabla}_X Y)\rho_{_x}.
\end{gathered}
\end{equation*}

Therefore
\begin{equation*}
\begin{aligned}
\nabla^2\rho_{_x}(X,Y) - \bar{\nabla}^2\rho_{_x} (X,Y)
=&\; (\bar{\nabla}_X Y - \nabla_X Y)\rho_{_x}\\
=&\; \vec{h}(X, Y)\rho_{_x}\\
=&\; \<\grad{M} \rho_{_x} \; ,\; \vec{h}(X, Y) \>.
\end{aligned}
\end{equation*}

Since $X$ and $Y$ are tangent to $L$, $Y\rho_{_x}$ is identically zero on $L$,
and so $XY\rho_{_x}(\Pi(x))=0$.  Therefore, at the point $\Pi(x)$,
\begin{equation*}
\begin{aligned}
\bar{\nabla}^2\rho_{_x}(X,Y)
=&\; X(Y\rho_{_x}) - (\bar{\nabla}_X Y)\rho_{_x}\\
=&\; -(\bar{\nabla}_X Y)\rho_{_x}\\
=&\; -\<\grad{M}\rho_{_x} \, ,\, \bar{\nabla}_X Y\>\\
=&\; -\norm{\grad{M} \rho_{_x}}\tilde{h}_{\Pi(x)}(X,Y).
\end{aligned}
\end{equation*}
It follows that, at the point $\Pi(x)$,
\begin{equation*}
\nabla^2\rho_{_x}(X,Y) = -\norm{\grad{M} \rho_{_x}}\tilde{h}_{\Pi(x)}(X,Y) + 
\<\grad{M} \rho_{_x}\, ,\, \vec{h}(X, Y) \>.
\end{equation*}

Since
$\norm{\grad{M} \rho_{_x}(\Pi(x))}=2d_M(x,N)$, it follows
from \cite[Theorem 6.4.3]{Petersen} that if $Y$ is a unit vector, then
\begin{equation*}
\begin{gathered}
(\nabla^2\rho_{_x})_{\Pi(x)}(Y,Y) \le
\frac{2d_M(x,N)\sqrt{\lambda}}{\tanh(d_M(x,N)\sqrt{\lambda})}+ 2\kappa d_M(x,N),
\quad \text{and}\\
(\nabla^2\rho_{_x})_{\Pi(x)}(Y,Y) \ge
\frac{2d_M(x,N)\sqrt{\lambda}}{\tan(d_M(x,N)\sqrt{\lambda})}- 2\kappa d_M(x,N).
\end{gathered}
\end{equation*}
\end{proof}

\begin{lem}\label{L:psi1}
There exists $\beta>0$ such that 
for all $g\in C^\infty(N)$,
\begin{equation*}  
\lim_{t \to 0} \frac{1}{(4\pi t)^{n/2}}
\int_N e^{\frac{-\rho_{_x}(y)}{4t}} g(y)\, d\nu(y)
= \frac{2^{n/2}}{\sqrt{\det(\nabla^2\rho_{_x})_{\Pi(x)}}} g(\Pi(x)),
\end{equation*}
uniformly on $N_\beta$.
\end{lem}

\begin{proof}
We use Laplace's method (see \cite[Theorem 15.2.2]{simon}).

Define $f:\R \to \R$ by
\begin{equation*}
f(u)=\frac{2\sqrt{\lambda}u}{\tan(\sqrt{\lambda}u)}- 2\kappa u, \quad u\in \R.
\end{equation*}
Then there exists $\beta'>0$ such that $f(u)\geq 1$ if $\abs{u}<\beta'$.
Put
\begin{equation*}
\beta=\min\left\{\frac{\eta}{2}, \beta'
\right\}.
\end{equation*}
Since $N$ is compact, there exists $C>0$ such that
\begin{equation}\label{E:third-covariant}
\abs{\nabla^3 \rho_{_x}(X,X,X)} \le C\norm{X}^3, \quad X\in TN.
\end{equation}

For $y\in N$ and $Y\in T_yN$, let
$J(y,Y)$ denote the absolute value of the Jacobian at $Y$ of
$\exp_y:T_yN \to N$.

Then $J:TN\to\R$ is smooth and $J(y,0)=1$ for all $y \in N$.

Since the function $g_1(y,Y)=g(\exp_y(Y))J(y,Y)$ is smooth, and $N$ is compact,
it follows that there exists a constant $B$ such that if
$(y,Y)\in TN$ then
\begin{equation*}
\begin{aligned}
\abs{g_1(y,Y)} \le B, \quad\text{and}\\
\abs{g_1(y,Y)-g(y)} \le B\norm{Y}.
\end{aligned}
\end{equation*}

Let $r=\min\{\inj{N}, 3/(2C)\}$.
Observe that the continuous function $\rho_{_x}(y)$ is strictly positive
on the compact set
\begin{equation*}
\{(x,y)\in M\x N \st
  d_M(x,N) \le \beta,
  d_N(\pi(x), y) \ge r\},
\end{equation*}
and hence there exists $b>0$ such that $\rho_{_x}(y) \ge b$ when
$x\in N_\beta$ and $d_N(\pi(x), y) \ge r$.
Therefore
\begin{equation*}
\lim_{t \to 0} \frac{1}{(4\pi t)^{n/2}}
\int_{d_N(\Pi(x),y) \ge r}
e^{\frac{-\rho_{_x}(y)}{4t}} g(y)\, d\nu(y) = 0,
\end{equation*}
uniformly on $N_\beta$.

Let $x\in N_\beta$ and suppose
$Y\in T_{\Pi(x)}N$ is a unit vector.
It follows from Lemma \ref{L:hessian} that
\begin{equation}\label{E:second-covariant}
\nabla^2\rho_{_x}(Y,Y) \ge f(d_M(x,N) \ge 1.
\end{equation}
Let $y\in N$
and suppose $d_N(\Pi(x),y) < r$. There exists 
a unique unit speed geodesic $c$ in $N$ such that $c(0)=\Pi(x)$
and $c(s)=y$ for some $s\in [0,r)$.
Let $Z= \dot{c}(0)$.
By Taylor's theorem, there exists $u \in (0,s)$ such that
\begin{equation*}
\begin{aligned}
\rho_{_x} (y)
  =&\; \(\rho_{_x} \circ c\)(s)\\
  =&\; \rho_{_x}(\Pi(x)) + s\(\nabla \rho_{_x}\)_{\Pi(x)}(Z) +
       \frac{s^2}{2}(\nabla^2 \rho_{_x})_{\Pi(x)}(Z,Z) +
       \frac{s^3}{6}{\nabla^3 \rho_{_x}(\dot{c}(u),\dot{c}(u),\dot{c}(u))}\\
  =&\; \frac{s^2}{2}(\nabla^2 \rho_{_x})_{\Pi(x)}(Z,Z) +
       \frac{s^3}{6}{\nabla^3 \rho_{_x}(\dot{c}(u),\dot{c}(u),\dot{c}(u))},
\end{aligned}
\end{equation*}
and hence, by equations (\ref{E:third-covariant}) and
(\ref{E:second-covariant}), it follows that
\begin{equation*}
\rho_{_x} (y) \ge s^2/4 = d_N(\Pi(x),y)^2/4.
\end{equation*}

First observe that
\begin{equation*}
\int_{T_{\Pi(x)}N} e^{-\frac{\nabla^2\rho_{_x}(Y,Y)}{8}} \,dY
 = \frac{(8\pi)^{n/2}}{\sqrt{\det(\nabla^2\rho_{_x})_{\Pi(x)}}}.
\end{equation*}

Let $B_r=\{Y \in T_{\Pi(x)}N \st \norm{Y} < r\}$.  Let $Y\in B_r$.

By the mean value theorem applied to the function $e^{-(\cdot)/4t}$,
and the fact that 
\begin{equation*}
\rho_{_x}(\exp_{\Pi(x)}(t^{1/2}Y))\ge \frac{t\norm{Y}^2}{4}, \quad\text{and}\quad
\frac{\nabla^2\rho_{_x}(t^{1/2}Y, t^{1/2}Y)}{2} \ge \frac{t\norm{Y}^2}{4},
\end{equation*}
we get
\begin{equation*}
\begin{aligned}
\abs{e^{-\frac{\rho_{_x}(\exp_{\Pi(x)}(t^{1/2}Y))}{4t}}
  - e^{-\frac{\nabla^2\rho_{_x}(t^{1/2}Y,t^{1/2}Y)/2}{4t}}}
\le&\;
\frac{e^{-\norm{Y}^2/16}}{4t}
\abs{
\rho_{_x}(\exp_{\Pi(x)}(t^{1/2}Y)) -
\nabla^2\rho_{_x}(t^{1/2}Y,t^{1/2}Y)/2
}\\
\le&\;
\frac{Ct^{1/2}\norm{Y}^3e^{-\norm{Y}^2/16}}{24}.
\end{aligned}
\end{equation*}

Then
\begin{equation*}
\begin{gathered}
  \abs{\frac{1}{(4\pi t)^{n/2}}\int_{d_N(\Pi(x),y) < r}
  e^{\frac{-\rho_{_x}(y)}{4t}} g(y)\, d\nu(y)
  - \frac{2^{n/2}g(\Pi(x))}{\sqrt{\det(\nabla^2\rho_{_x})_{\Pi(x)}}}}\\
= \abs{\frac{1}{(4\pi t)^{n/2}}
  \int_{B_r} e^{-\frac{\rho_{_x}(\exp_{\Pi(x)}(Y))}{4t}} g_1(\Pi(x),Y) \,dY
  - \frac{2^{n/2}g(\Pi(x))}{\sqrt{\det(\nabla^2\rho_{_x})_{\Pi(x)}}}}\\
= \abs{\frac{1}{(4\pi)^{n/2}}
  \int_{B_{t^{-1/2}r}}
  e^{-\frac{\rho_{_x}(\exp_{\Pi(x)}(t^{1/2}Y))}{4t}} g_1(\Pi(x),t^{1/2}Y) \,dY
  - \frac{2^{n/2}g(\Pi(x))}{\sqrt{\det(\nabla^2\rho_{_x})_{\Pi(x)}}}}\\
= \frac{1}{(4\pi)^{n/2}} \abs{
  \begin{gathered}
  \int_{B_{t^{-1/2}r}}
  \(
  e^{-\frac{\rho_{_x}(\exp_{\Pi(x)}(t^{1/2}Y))}{4t}} g_1(\Pi(x),t^{1/2}Y)
  - e^{-\frac{\nabla^2\rho_{_x}(t^{1/2}Y,t^{1/2}Y)/2}{4t}} g(\Pi(x)) 
  \)  \, dY\\
  -  g(\Pi(x)) \int_{B_{t^{-1/2}r}^c} e^{-\frac{\nabla^2\rho_{_x}(Y,Y)}{8}} \, dY
  \end{gathered}
}\\
= \frac{1}{(4\pi)^{n/2}} \abs{
  \begin{gathered}
  \int_{B_{t^{-1/2}r}}
  e^{-\frac{\rho_{_x}(\exp_{\Pi(x)}(t^{1/2}Y))}{4t}}
  \(g_1(\Pi(x),t^{1/2}Y) - g(\Pi(x))\) \,dY\\
  +
  \int_{B_{t^{-1/2}r}}
  \(e^{-\frac{\rho_{_x}(\exp_{\Pi(x)}(t^{1/2}Y))}{4t}}
  - e^{-\frac{\nabla^2\rho_{_x}(t^{1/2}Y,t^{1/2}Y)/2}{4t}}\)
 g(\Pi(x)) \,dY\\
  -  g(\Pi(x)) \int_{B_{t^{-1/2}r}^c} e^{-\frac{\nabla^2\rho_{_x}(Y,Y)}{8}} \, dY
  \end{gathered}
  }\\
\le
\frac{2^{(n+2)}\Gamma((n+1)/2)\(3B+C(n+1)g(\Pi(x))\)}{3\Gamma(n/2)} t^{1/2}
+ \frac{2^{n/2} g(\Pi(x))}{\sqrt{\det(\nabla^2\rho_{_x})_{\Pi(x)}}} e^{-nr^2/16t}.
\end{gathered}
\end{equation*}
Therefore
\begin{equation*}
\lim_{t \to 0} \frac{1}{(4\pi t)^{n/2}}\int_{d_N(\Pi(x),y) < r}
            e^{\frac{-\rho_{_x}(y)}{4t}} g(y)\, d\nu(y)
= \frac{2^{n/2}g(\Pi(x))}{\sqrt{\det(\nabla^2\rho_{_x})_{\Pi(x)}}},
\end{equation*}
uniformly on $N_\beta$.
\end{proof}

\section{The heat flow of $\tau$}\label{S:heat-flow}

\begin{thm}\label{T:heat}
Let $f_t$ be the heat flow of $\tau$. Then
\begin{equation*}
\norm{f_{t/2}}_2^2 \sim (4\pi t)^{-k/2}\norm{\psi}_2^2, \quad t\to 0.
\end{equation*}
\end{thm}

The Landau notation will be used throughout, and will always refer to
the asymptotic as $t\to 0^+$.
By Minakshisundaram and Pleijel's asymptotic formula for the heat kernel
(see \cite[\S1]{minakshisundaram}, \cite[\S4]{patodi} or
\cite[Theorem 7.15]{roe}), there
exists a number $l$ and smooth functions $\Theta_i \in C^\infty(M\x M)$,
$i=1,\dots, l$, with $\Theta_0(x,x)=1$ such that
\begin{equation*}
k_t(x,y) = \frac{e^{-d_M(x,y)^2/4t}}{{(4\pi t)}^{m/2}}
(\Theta_0(x,y) + t\Theta_1(x,y) + t^2\Theta_2(x,y) + \cdots + t^l\Theta_l(x,y))
+ O(t),
\end{equation*}
uniformly on $M\x M$.
Therefore
\begin{equation*}
k_t(x,y) = \frac{e^{-d_M(x,y)^2/4t}}{{(4\pi t)}^{m/2}}\Theta_0(x,y)
+ o\(1 + t^{-m/2}e^{\frac{-d_M(x,y)^2}{4t}}\),
\end{equation*}
uniformly on $M\x M$.

Let
\begin{equation*}
f_t^0(x)
= \frac{1}{(4\pi t)^{m/2}}
  \int_N e^{-\frac{d_M(x,y)^2}{4t}}\, \Theta_0(x,y) \psi(y) \,d\nu(y),
  \quad x \in M.
\end{equation*}

Observe that
\begin{equation*}
f_t(x) = f_t^0(x) + o\(\int_N t^{-m/2} e^{-\frac{d_M(x,y)^2}{4t}} \psi(y) d\nu(y)\),
\end{equation*}
uniformly on $M$.

By Lemma \ref{L:psi1},
\begin{equation*}
\begin{gathered}
f_t(x) = f_t^0(x) + o\(t^{-k/2}e^{-\frac{d_M(x,N)^2}{4t}}\),  \quad\text{and}\\
f_t^0(x) = \frac{2^{n/2} \Theta_0(x,\Pi(x))\psi(\Pi(x))e^{-d_M(x,N)^2/4t}}
                {(4\pi t)^{k/2}\sqrt{\det(\nabla^2\rho_{_x})_{\Pi(x)}}}
           + o\(t^{-k/2}e^{-\frac{d_M(x,N)^2}{4t}}\),
\end{gathered}
\end{equation*}
uniformly on $N_\beta$.
Let $\boldsymbol{\nu}$ denote the normal bundle of $N$ in $M$, $dV$
the volume element in $\boldsymbol{\nu}$,
$\exp_{\boldsymbol{\nu}} : \boldsymbol{\nu} \to M$ the corresponding
exponential map, and $\vartheta_v(s)$ the infinitesimal change of volume
function of $N$ in the direction of the unit vector $v$ (see
\cite[\S3.2]{tubes}).  By \cite [33.2.2]{Burago},
\begin{equation*}
dV(y,sv) = s^{k-1} \, d\nu(y) \, d\sigma(v) \, ds,
\end{equation*}
where $d\sigma$ is the volume element of the unit sphere $\Sigma_{k-1}$.
Therefore
\begin{equation*}
\begin{gathered}
\norm{\frac{2^{n/2}e^{-d_M(\cdot,N)^2/4t} \Theta_0(\cdot,\Pi(\cdot)) \psi(\Pi(\cdot))}
  {(4\pi t)^{k/2}\sqrt{\det{(\nabla^2\rho_{_{(\cdot)}})_{\Pi(\cdot)}}}}}_{L^2(N_\beta)}^2\\
=\frac{2^n}{(4\pi t)^k}
\int_{\boldsymbol{\nu}}
\chi_{\exp_{\boldsymbol{\nu}}^{-1}(N_\beta)}
\frac{e^{-s^2/2t} \Theta_0(\exp_{\boldsymbol{\nu}}(y, sv), y)^2 \abs{\psi(y)}^2}
     {\det{(\nabla^2\rho_{_{\exp_{\boldsymbol{\nu}}(y, sv)}})_y}}
     \vartheta_v(s)\, dV(y,sv)\\
=\frac{2^n}{(4\pi t)^k}
\int_N \int_0^\beta \int_{\Sigma_{k-1}}
\frac{e^{-s^2/2t} \Theta_0(\exp_{\boldsymbol{\nu}}(y, sv), y)^2 \abs{\psi(y)}^2}
     {\det{(\nabla^2\rho_{_{\exp_{\boldsymbol{\nu}}(y, sv)}})_y}} \vartheta_v(s)
     s^{k-1} \, d\sigma(v) \, ds \, d\nu(y)\\
=\frac{2^n}{(4\pi t)^k}
\int_N \int_0^{t^{-1/2}\beta} \int_{\Sigma_{k-1}}
\frac{e^{-s^2/2} \Theta_0(\exp_{\boldsymbol{\nu}}(y,t^{1/2}sv), y)^2 \abs{\psi(y)}^2}
     {\det{(\nabla^2\rho_{_{\exp_{\boldsymbol{\nu}}(y, t^{1/2}sv)}})_y}}
     \vartheta_v(t^{1/2}s) t^{k/2} s^{k-1} \, d\sigma(v) \, ds \, d\nu(y).
\end{gathered}
\end{equation*}
Therefore, by Lemma \ref{L:hessian} and
\cite[Lemma 8.24 and Corollary 8.26]{tubes}, we have
\begin{equation*}
\norm{\frac{2^{n/2}e^{-d_M(\cdot,N)^2/4t} \Theta_0(\cdot,\Pi(\cdot)) \psi(\Pi(\cdot))}
  {(4\pi t)^{k/2}\sqrt{\det{(\nabla^2\rho_{_{(\cdot)}})_{\Pi(\cdot)}}}}}_{L^2(N_\beta)}^2
=
\frac{\norm{\psi}_2^2}{(8\pi t)^{k/2}} + o(t^{-k/2}).
\end{equation*}
Since
$
\int_{M\setminus N_\beta} \abs{f_t^0}^2 \,d\mu = o(1),
$
it follows that
\begin{equation*}
\norm{f_t^0}_2 = 
\frac{\norm{\psi}_2}{(8\pi t)^{k/4}}
+ o(t^{-k/4}).
\end{equation*}
Since
$
\int_{M\setminus N_\beta} \abs{f_t}^2 \,d\mu = o(1),
$
it follows that
\begin{equation*}
\begin{aligned}
\abs{\norm{f_t}_2-\frac{\norm{\psi}_2}{(8\pi t)^{k/4}}}
  =&\; \abs{\norm{f_t}_2 - \norm{f_t^0}_2} + o(t^{-k/4})\\
\le&\; \norm{f_t-f_t^0}_2 + o(t^{-k/4})\\
\le&\; o\(t^{-k/4}\).
\end{aligned}
\end{equation*}
Therefore
\begin{equation*}
\norm{f_t}_2 \sim \frac{\norm{\psi}_2}{(8\pi t)^{k/4}}, \quad t\to 0.
\end{equation*}
Therefore
\begin{equation*}
\norm{f_{t/2}}_2^2 \sim \frac{\norm{\psi}_2^2}{(4\pi t)^{k/2}}, \quad t\to 0.
\end{equation*}

\section{Asymptotic Pythagorean identity for the Legendre Polynomials}\label{S:lp}

Suppose $M$ is the unit sphere in $\R^3$, and $N$ is a latitude
circle,
say
$$
\begin{aligned}
M=&\; \{(x,y,z)\in\R^3 \st x^2+y^2+z^2=1\}, \quad\text{and}\\
N=&\; \{(x,y,z)\in M \st z=z_0\}.
\end{aligned}
$$

Let $\tau=\nu$ be the length measure on $N$ (i.e. $\psi=1$).
Then, according to Theorem \ref{T:main},
\begin{equation}\label{E:sph-harm}
\begin{aligned}
\sum_{\lambda_j < \lambda} \abs{\hat{\tau}(j)}^2
   =&\; \sum_{\lambda_j < \lambda} \abs{\<\tau, \phi_j\>}^2
\sim&\; \frac{\nu(N)}{\pi} \lambda^{1/2} 
   =&\; 2\sqrt{1-z_0^2} \lambda^{1/2}.
\end{aligned}
\end{equation}

Now suppose the $\{\phi_j\}$ are the standard spherical harmonics
$\{\{Y_l^m\}_{m=-l}^{l}\}_{l=0}^\infty$.
Since $\tau$ is zonal, $\<\tau, Y_l^m\>=0$ if $m\ne 0$.
Recall that
$$
-\Delta Y_l^m = l(l+1) Y_l^m.
$$
Therefore (\ref{E:sph-harm}) reads
\begin{equation}\label{E:sph-harm-2}
\sum_{l(l+1)<\lambda} \abs{\<\tau, Y_l^0\>}^2 \sim 2\sqrt{1-z_0^2} \lambda^{1/2}.
\end{equation}

Since $Y_l^0(x,y,z)=\sqrt{\frac{2l+1}{4\pi}} P_l(z)$, where $P_l(z)$ is the
$l$-th Legendre polynomial,
Equation (\ref{E:sph-harm-2}) gives
\begin{equation*}
\sum_{l(l+1)<\lambda} \(l+\frac12\)\abs{P_l(z_0)}^2 \sim
\frac{\lambda^{1/2}}{\pi\sqrt{1-z_0^2}}
\end{equation*}
or equivalently
\begin{equation*}
\lim_{N\to\infty} \frac{1}{N} \sum_{l=0}^N \(l+\frac12\) \abs{P_l(z_0)}^2
= \frac{1}{\pi\sqrt{1-z_0^2}},
\end{equation*}
which is an asymptotic Pythagorean identity for the Legendre
Polynomials.
It can also be proved quite easily using the
Christoffel-Darboux summation formula (see \cite[Equation 3.2.4]{szego})
and Laplace's asymptotic formula for the Legendre polynomials 
(see \cite[Theorem 8.21.2]{szego}).

Analogous formulae for the Jacobi polynomials $P_n^{(\alpha,\beta)}$
(for $\alpha,\beta$ appearing in Table 1 of \cite{ali}) may
be derived by interpreting them as spherical functions of
rank-one symmetric spaces of compact type.

\bibliographystyle{amsplain}
\bibliography{v14-grmp}

\end{document}